\documentclass[a4paper,11pt]{article}
\usepackage{latexsym}
\usepackage{amsmath}
\usepackage{amssymb}
\usepackage{enumerate}
\usepackage{theorem}
\usepackage{array}
\newtheorem{theorem}{Theorem}

\theorembodyfont{\rmfamily}

\newtheorem{setting}[theorem]{Setting}

\newcommand{\proof}{\noindent \mbox{\em Proof.\hspace*{2mm}}}
\newcommand{\qed}{\hfill \mbox{$  \Box $}}

\DeclareMathOperator{\image}{Im}

\title{A bound for the number of lines lying on
a nonsingular surface in $3$-space
over a finite field}
\author{
Masaaki Homma
\thanks{Partially supported by Grant-in-Aid
for Scientific Research (24540056, 15K04829), Japan Society for the Promotion of Science (JSPS).}
\\
 Department of Mathematics and Physics\\
Kanagawa University\\
Hiratsuka 259-1293, Japan\\
homma@kanagawa-u.ac.jp
\and
Seon Jeong Kim
\thanks{Partially supported by Basic Science Research Program through the National Research Foundation of Korea (NRF) 
funded by the Ministry of Education, Science and Technology (2012R1A1A2042228).
}\\
 Department of Mathematics and RINS\\
Gyeongsang National University\\
Jinju 660-701, Korea \\
skim@gnu.kr
}
\date{}

\begin{document}
\maketitle
\begin{abstract}
A nonsingular surface of degree $d \geq 2$ in $\mathbb{P}^3$
over $\mathbb{F}_q$ has at most
$((d-1)q+1)d$ $\mathbb{F}_q$-lines,
and this bound is optimal for $d = 2, \sqrt{q}+1, q+1$.
This is a bi-product of a previous study on estimating
the number of $\mathbb{F}_q$-points of surfaces.
\end{abstract}

This short note is a supplement to our previous works
\cite{hom-kim2013b, hom-kim2015online, hom-kim2015b}.
In \cite{hom-kim2015online, hom-kim2015b},
we considered the set of $\mathbb{F}_q$-lines on a given surface
$S \subset \mathbb{P}^3$
as an auxiliary tool in order to estimate the number of
$\mathbb{F}_q$-points of $S$.
This time, we turn our attention to $\mathbb{F}_q$-lines themselves
on $S$.

\begin{setting}\label{setting}
Let $S$ be a nonsingular surface of degree $d \geq 2$ in $\mathbb{P}^3$
defined over $\mathbb{F}_q$.
Let $\nu_q(S)$ denote the number of $\mathbb{F}_q$-lines lying on $S$,
and $N_q(S)$ that of $\mathbb{F}_q$-points of $S$.
\end{setting}

In \cite{hom-kim2013b, hom-kim2015b},
we showed the following fact, in which $N_q(X)$ denotes 
the number of $\mathbb{F}_q$-poins of the surface $X$,
though $X$ is not necessary
nonsingular.

\begin{theorem}\label{theoremonpoints}
Let $X$ be a surface in ${\Bbb P}^3$ over ${\Bbb F}_q$
without ${\Bbb F}_q$-plane components,
then
\begin{equation}\label{elementarybound}
 N_q(X) \leq ((d-1)q +1)(q+1).
\end{equation}
Furthermore, equality holds in {\rm (\ref{elementarybound})},
if and only if 
the surface $X$ is projectively equivalent to
one of the following surfaces over ${\Bbb F}_q${\rm :}
\begin{enumerate}[{\rm (i)}]
\item
$X_0X_1-X_2X_3=0$  if $d=2${\rm ;}
\item
$X_0^{\sqrt{q}+1}+ X_1^{\sqrt{q}+1} +X_2^{\sqrt{q}+1}+X_3^{\sqrt{q}+1}=0$
if $d= \sqrt{q}+1${\rm ;}
\item $X_0X_1^{q} - X_0^{q}X_1 + X_2X_3^{q} - X_2^{q}X_3 =0$
if $d=q+1$.
\end{enumerate}
\end{theorem}
In this note, we shall give a proof of the following theorem.
In the proof, ${}^{\#}T$ denotes the cardinality of $T$
if $T$ is a finite set.
\begin{theorem}\label{maintheorem}
Under Setting~{\rm \ref{setting}},
\begin{equation}\label{boundforlines}
\nu_q(S) \leq \frac{d}{q+1} N_q(S) \leq ((d-1)q +1)d.
\end{equation}
Furthermore the list of surfaces $S$ satisfying
$\nu_q(S) = ((d-1)q +1)d$ coincides with
that in Theorem~{\rm \ref{theoremonpoints}}.
\end{theorem}
\proof
Let $G(1, \mathbb{P}^3)$ be the Grassmann variety of
lines in $\mathbb{P}^3$.
The sets of $\mathbb{F}_q$-points of $S$ and $G(1, \mathbb{P}^3)$
are denoted by $S(\mathbb{F}_q)$ and $G(1, \mathbb{P}^3)(\mathbb{F}_q)$
respectively.
Consider the correspondence
\[
\Pi = \{ (P, l) \mid P \in l \subset S \}
 \subset S(\mathbb{F}_q) \times G(1, \mathbb{P}^3)(\mathbb{F}_q)
\]
with projections
$\pi_1 : \Pi \to S(\mathbb{F}_q)$
and
$\pi_2: \Pi \to G(1, \mathbb{P}^3)(\mathbb{F}_q)$.
Then
${}^{\#}( \image \pi_2) = \nu_q(S)$,
and hence ${}^{\#} \Pi = (q+1) \nu_q(S).$
On the other hand,
a line $l$ on $S$ passing through $P$ lies on the tangent plane $T_P(S)$
to $S$ at $P$ because $l = T_P(l) \subset T_P(S)$.
Hence the number of lines on $S$ passing through the assigned point $P$
is at most that of line components of the curve $S \cap T_P(S)$.
Hence
${}^{\#} \pi_1^{-1}(P) \leq d$.
So, to sum up,
$(q+1) \nu_q(S) ={}^{\#} \Pi \leq N_q(S)d,$
which is the first inequality in (\ref{boundforlines}).
The second one comes from Theorem~\ref{theoremonpoints}.

Next we show the additional statement.
If $\nu_q(S) = ((d-1)q +1)d$,
then $N_q(S) = ((d-1)q +1)(q+1)$
by (\ref{boundforlines}),
and hence $S$ is one of the surfaces listed in Theorem~\ref{theoremonpoints}.
Conversely we show that the following claim holds for
each surface $S$ defined by (i) or (ii) or (iii);\\
Claim: for any $P \in S(\mathbb{F}_q)$, the intersection of
the tangent plane $T_P(S)$ with $S$ is a union of $d$
$\mathbb{F}_q$-lines with vertex $P$, that is,
$T_P(S) \cap S$ forms a planar pencil with vertex $P$
in terms of \cite{hom-kim2015b}.
This claim actually implies
$\nu_q(S) = ((d-1)q + 1)d$.
Indeed,
${}^{\#} \pi_1^{-1}(P) = d$
in the first part of this proof,
together with
$N_q(S) = ((d-1)q +1)(q+1)$,
gives rise to
\[
 (q+1) \nu_q(S)={}^{\#} \Pi = N_q(S)d =  ((d-1)q + 1)(q+1)d.
\]
This claim had been proved essentially in \cite[Prop.~3.1]{hom-kim2015b}, however, here we give a direct proof by using equations (i), (ii) or (iii).
Let $P=(a_0, \dots , a_3) \in S(\mathbb{F}_q)$.

(i) Suppose that $S$ is defined by equation (i).
Then the tangent plane $T_P(S)$
is defined by
\[
a_1X_0 +a_0 X_1 - a_3 X_2 -a_2 X_3 =0.
\]
Without loss of generality, we may assume that $a_1=1$, hence
$a_0=a_2a_3$.
Then $S(\mathbb{F}_q) \cap T_P(S)$
is defined by
\begin{align*}
0 &= (-a_2a_3 X_1 + a_3X_2 + a_2 X_3)X_1 -X_2 X_3 \\
  &= (a_2 X_1-X_2)(X_3 -a_3 X_1)
\end{align*}
in $T_P(S) =\mathbb{P}^2$ with coordinates $X_1, X_2, X_3$.

(ii) Suppose that $S$ is defined by equation (ii).
Since any automorphism of $S$ comes from an $\mathbb{F}_q$-linear
transformation of $\mathbb{P}^2$
and the automorphism group acts on $S(\mathbb{F}_q)$ as transitively
\cite[Lemmas 3.7 and 3.8]{hom2006},
we may assume that
$P=(1, \zeta , 0, 0)$
with $\zeta^{\sqrt{q}+1}=-1$.
Hence $T_P(S)$ is defined by $X_0 + \zeta^{\sqrt{q}}X_1=0$,
and $T_P(S) \cap S$ by $X_2^{\sqrt{q}+1} + X_3^{\sqrt{q}+1}=0$
in $T_P(S) =\mathbb{P}^2$ with coordinates $X_1, X_2, X_3$.
Since 
\[X_2^{\sqrt{q}+1} + X_3^{\sqrt{q}+1}
= \prod_{
\begin{subarray}{c}
\lambda \in \mathbb{F}_q^{\ast}\\
\text{with} Nm\, \lambda = 1
\end{subarray}
}
(X_2 - \lambda \zeta X_3),
\]
$T_P(S) \cap S$ splits into $\sqrt{q}+1$ lines.

(iii) Suppose that $S$ is defined by equation (iii).
Since $S(\mathbb{F}_q)=\mathbb{P}^3(\mathbb{F}_q)$,
it is enough to show that
for $Q=(b_0, \dots, b_3) \in S(\mathbb{F}_q) \cap T_P(S)$,
the line 
$PQ=\{ \lambda(a_0, \dots , a_3)  + \mu (b_0, \dots, b_3)
     \mid (\lambda , \mu ) \in \mathbb{P}^1 \}$
joining $P$ and $Q$ lies on $S$.
Since $P=(a_0, \dots , a_3)$ is
an $\mathbb{F}_q$-point, the tangent plane $T_P(S)$
is defined by
\[
a_1X_0 - a_0 X_1 + a_3 X_2 -a_2 X_3 =0.
\]
Hence
\[
a_1b_0 - a_0b_1 + a_3 b_2 -a_2 b_3 =0.
\]
Therefore
\begin{multline*}
(\lambda a_0 + \mu b_0)(\lambda a_1 + \mu b_1)^q
- (\lambda a_0 + \mu b_0)^q(\lambda a_1 + \mu b_1)\\
\shoveright{+
(\lambda a_2 + \mu b_2)(\lambda a_3 + \mu b_3)^q
- (\lambda a_2 + \mu b_2)^q(\lambda a_3 + \mu b_3)}\\
\shoveleft{
= \lambda^q \mu (a_1b_0 - a_0b_1 + a_3 b_2 -a_2 b_3)
 + \lambda \mu^q (a_0b_1 - a_1b_0 + a_2 b_3 -a_3 b_2)}
\end{multline*}
is identically $0$.
\qed


\begin{thebibliography}{00}
\bibitem{hom2006}
M. Homma,
{\em Galois points for a Hermitian curve},
Comm. Algebra 34 (2006) 4503--4511.
\bibitem{hom-kim2013b}
M. Homma and S. J. Kim,
{\em An elementary bound for the number of points of
a hypersurface over a finite field},
Finite Fields Appl. 20 (2013) 76--83.
\bibitem{hom-kim2015online}
M. Homma and S. J. Kim,
{\em The characterization of Hermitian surfaces by the number of points},
J. Geom. (on line) 18 July 2015, DOI 10.1007/s00022-015-0283-1.
\bibitem{hom-kim2015b}
M. Homma and S. J. Kim,
{\em Numbers of points of surfaces in the projective $3$-space
over finite fields},
Finite Fields Appl. 35 (2015) 52--60.

\end{thebibliography}
\end{document}